\documentstyle{ajour}

\makeatletter
   \@addtoreset{equation}{section}
   \renewcommand{\theequation}{\thesection.\arabic{equation}}
   \let\savetheequation\theequation
   \makeatother

\begin{document}

\title{On the questions P $\stackrel{?}{=}$ NP $\cap$ co-NP and NP $\stackrel{?}{=}$ co-NP for infinite time Turing machines}

\titlerunninghead{P $=$ NP $\cap$ co-NP for infinite time Turing machines}
\author{Vinay Deolalikar}
\authorrunninghead{Vinay Deolalikar}
\affil{Hewlett-Packard Research \\ 1501 Page Mill Road, M/S 3U-4
\\ Palo Alto, CA 94304 \\
Tel: (650) 857 8605 Fax: (650) 852 3791}
\email{vinayd@hpl.hp.com}

\begin{article}

\renewcommand{\labelenumi}{$(\roman{enumi})\;\;$}

\newtheorem{mydefinition}[theorem]{Definition}
\newtheorem{mylemma}[theorem]{Lemma}
\newtheorem{myproposition}[theorem]{Proposition}
\newtheorem{mycorollary}[theorem]{Corollary}
\newtheorem{myexample}[theorem]{\sc Example}
\newtheorem{myremark}[theorem]{Remark}

\abstract{Schindler recently addressed two versions of the
question P $\stackrel{?}{=}$ NP for Turing machines running in
transfinite ordinal time. These versions differ in their
definition of input length. The corresponding complexity classes
are labelled P, NP and ${\rm P}^+,{\rm NP}^+$. Schindler showed
that P $\neq$ NP and ${\rm P}^+ \neq {\rm NP}^+$. We show that P
$=$ NP $\cap$ co-NP and NP $\neq$ co-NP, whereas ${\rm P}^+
\subset$ NP $\cap$ co-NP and ${\rm NP}^+ \neq $ co-NP$^+$.
}

\keywords{Complexity theory, Descriptive set theory.}

\section{Introduction} \label{section:introduction}

The fundamental open problems in complexity theory for Turing
machines running in finite time are whether P $\stackrel{?}{=}$
NP, whether P $\stackrel{?}{=}$ NP $\cap$ co-NP and whether NP
$\stackrel{?}{=}$ co-NP.

After Hamkins and Lewis \cite{HamLew} formalized the notion of a
Turing machine running in transfinite ordinal time, the
corresponding questions about infinite time Turing machines were
posed. Schindler recently defined two versions of the classes P
and NP in a Turing machine running in transfinite ordinal time.
He labelled them P, NP and ${\rm P}^+,{\rm NP}^+$ and showed
that P $\neq$ NP  and ${\rm P}^+ \neq {\rm NP}^+$. The
proofs used a classical theorem from descriptive set theory
stating that not all analytic sets are Borel. Analogies between
the classes P and NP, and the classes of Borel and analytic sets 
respectively, had earlier been drawn by Sipser and others 
(cf. \cite{Sip1} and the references therein).

We address the problems of whether P $\stackrel{?}{=}$ NP $\cap$
co-NP (resp. ${\rm P}^+ \stackrel{?}{=} {\rm NP}^+ \cap $
co-NP$^+$) and whether NP $\stackrel{?}{=}$ co-NP (resp. $ {\rm NP
}^+ \stackrel{?}{=} $ co-NP$^+$). We use a classical result by
Suslin which characterizes analytic sets that are Borel as those whose
complements are also analytic. Using Suslin's theorem, we prove
that P $=$ NP $\cap$ co-NP. We then show that ${\rm P}^+$ is a
strict subset of ${\rm NP}^+ \cap $ co-NP$^+$ using properties of
projections of hyperarithmetical sets.

We also observe that NP $\neq$ co-NP and $ {\rm NP }^+ \neq$ {\rm
co-NP}$^+$.

\section{Preliminaries}

This section will provide some notions from descriptive set theory
needed in order to state our results. For details, the reader is
referred to \cite{Mos,Kec}.

We first fix our notation. The first infinite ordinal is denoted
by $\omega$, the first uncountable ordinal by $\omega_1$, the
first non-recursive ordinal ({\it i.e.} the Church-Kleene ordinal)
by $\omega_1^{CK}$, and the ordinal of G\"{o}del's constructible
universe $L$ by $\omega_1^L$. A Polish space, {\it i.e.} a space
that is separable and completely metrizable, is denoted by $X$. In
this paper, we are particularly interested in Cantor space
$^{\omega}2$ and Baire space $^{\omega}\omega$.

The class of Borel sets in $X$ are denoted by ${\bf B}(X)$. This
class ramifies in an infinite Borel Hierarchy whose classes are
denoted by ${\bf \Sigma}_\xi^0$ and ${\bf \Pi}_\xi^0$. The
projective sets obtained from Borel sets by the operations of
projection and complementation ramify into a projective hierarchy
of length $\omega$ whose classes are denoted by ${\bf
\Sigma}_{n}^1$ and ${\bf \Pi}_{n}^1$. The ambiguous classes are
denoted by ${\bf \Delta}_{n}^1 = {\bf \Sigma}_n^1 \cap {\bf
\Pi}_n^1$. Members of the class ${\bf \Sigma}_1^1$ are called
analytic, members of ${\bf \Pi}_{1}^1$ are called co-analytic and
members of ${\bf \Delta}_{1}^1$ are called bi-analytic. The
effective analogues of the Borel classes are denoted by lightface
${\Sigma}_{\xi}^0$ and ${\Pi}_{\xi}^0$. Similarly the effective
analogues of the projective classes are denoted by
${\Sigma}_{n}^1$, ${\Pi}_{n}^1$ and ${\Delta}_{n}^1$.

We need the following two central results in descriptive set
theory.
\begin{theorem}[Lusin's Separation Theorem]
If $X$ is a standard Borel space, and $A,B \subseteq X$ are two
disjoint analytic sets, then there is a Borel set $C \subseteq X$
such that $A \subseteq C$ and $C \cap B = \emptyset$.
\end{theorem}

\begin{theorem}[Suslin's Theorem] For a Polish space $X$,
{\bf B}(X) = ${\bf \Delta}_1^1(X)$.
\end{theorem}
\begin{proof} The Effros Borel space of the set of closed subsets of a
Polish space is standard so the Lusin separation theorem can be
applied. The proof of Suslin's theorem is immediate by letting $B$ be the complement of $A$ in the Lusin separation theorem.
\end{proof}

\begin{definition}
Let $\alpha$ be a countable ordinal. We say that $A \subset
^{\omega}2$ is $\Delta^1_1(\alpha)$ in $\alpha$ if it is uniformly
$\Delta^1_1$ in any real $x$ coding a well-order of order type
$\alpha$.
\end{definition}

\begin{lemma} \label{lemma:proj_is_countable}
Let $A$ be $\Delta^1_1(\alpha)$ in a countable ordinal $\alpha$ as defined
above. Let $B$ be its projection. If $B$ is a Borel set, then there exists
a countable ordinal $\beta$ such that $B$ is $\Delta^1_1(\beta)$.
\end{lemma}
\begin{proof}
Let $WO$ be the set of codes of countable ordinals, and
$WO(\alpha)$ be the set of codes of the ordinal $\alpha$. Let $B^c$
be the complement of $B$. Since $B$ is $\Sigma^1_1$ in $\alpha$,
it follows that there is a $\Pi^1_1$ set $D$ such that
$$ x \in B^c \iff (x,w) \in D,  \mbox{ for any real } w \mbox{ coding } \alpha.$$
Let $f$ be a recursive function such that
$$(x,w) \in P \iff f(x,w) \in WO.$$
Then $B^c \times WO(\alpha)$ is Borel and contained in $D$,
therefore its image under $f$ is analytic and contained in $WO$.
By the Boundedness theorem, there is an ordinal $\delta$, such
that if $|w|$ denotes the ordinal coded by $w$, we have
$$x \in B^c \iff |f(x,w)| < |v|$$
for any code $w$ of $\alpha$ and any code $v$ of $\delta$. It immediately follows that $B$ is $\Delta^1_1(\beta)$, for some $\beta.$
\end{proof}

The countable ordinal $\beta$ is possibly much larger than
$\alpha$. We have the following result from \cite[Theorems 1.4,2.3]{KMS}.

\begin{lemma} \label{lemma:notless}
There exist $\Delta^1_1$ sets with projections that are
$\Delta_1^1(\omega_1^L)$, but not $\Delta^1_1(\alpha)$ for any
$\alpha < \omega_1^L$.
\end{lemma}

\begin{corollary} \label{corollary:notless}
There exist hyperarithmetical sets whose projections are Borel but
not hyperarithmetical.
\end{corollary}
\begin{proof} Immediate from the Lemma~\ref{lemma:notless} by letting $\alpha$ be recursive and observing that  $\omega_1^L > \omega_1^{CK}$.
\end{proof}

\section{Infinite time Turing machines}

Infinite time Turing machines were introduced in \cite{HamLew}.
There has been growing research interest in these machines
especially after Schindler \cite{Sch} showed that P $\neq$ NP in
this transfinite setting. We recall the following definitions from
\cite{Sch}.

\begin{definition}
Let $A \subset\,^{\omega}2$ and let $\alpha \leq \omega_1 + 1$. Then $A$ is in $P_\alpha$ if there exists a Turing machine $T$ and some $\beta < \alpha$ such that\\
(a) T decides A \\
(b) T halts on all inputs after $< \beta$ many steps.
\end{definition}

\begin{definition}Let $A \subset\, ^{\omega}2$ and let $\alpha \leq \omega_1 + 1$. Then $A$ is in $NP_\alpha$ if there exists a Turing machine $T$ and some $\beta < \alpha$ such that \\
(a) $x \in A$ if and only if $(\exists y$ such that $T$ accepts $x \oplus y)$ \\
(b) $T$ halts on all inputs after $< \beta$ many steps.
\end{definition}

\subsection{The classes P and NP in infinite time Turing machines}

If we let all inputs $x \in\, ^{\omega}2$ as having the same
length $\omega$, then we have the following characterization of
the classes P and NP.
\begin{definition} P $= P_{\omega^\omega }$ and NP $= NP_{\omega^\omega }$.
\end{definition}

The following description of the class $P$ is given in
\cite[Lemmas 2.5, 2.6]{Sch}.
\begin{lemma} Let $A \subset\, ^{\omega}2$.
Then $A \in P_{\omega_1^{CK}}$ if and only if $A$ is a
hyperarithmetic set. Furthermore, $A \in P_{\omega_1}$ if and only
if A is $\Delta_1^1(\alpha)$ in a countable ordinal $\alpha$.
\end{lemma}

In particular, for a Borel set to be outside of $P_{\omega_1}$, it
must be $\Delta_1^1$ in a real $x$ that does not code any
countable well-order.

Now we state our main theorem.
\begin{theorem} \label{theorem:NPintersectcoNP}
{\rm P = NP $\cap$ co-NP.}
\end{theorem}
\begin{proof}
Since P $\subset$ NP and P $=$ co-P, it follows that P $\subset$
co-NP and therefore P $\subseteq$ NP $\cap$ co-NP. Let $B \subset\,
^{\omega}2$ satisfy $B \in$ NP $\cap$ co-NP. Since $B \in$ NP, it 
is the projection of a set in P. However, P $\subset {\bf \Delta}_1^1$, 
which implies that $B \in {\bf \Sigma}_1^1$. Let $B^{c}$ be the complement of
$B$. Since $B \in$ co-NP, it follows that $B^{c} \in NP$. Thus
$B^{c} \in {\bf \Sigma}_1^1$, which implies that $B \in {\bf
\Pi}_1^1$. Therefore $B \in  {\bf \Sigma}_1^1 \cap {\bf \Pi}_1^1 =
{\bf \Delta}_1^1$. From Suslin's theorem, it follows that $B$ must
be Borel.

We cannot immediately infer that $B$ is in P because not all Borel
sets are in P. Those in P are precisely those that are
$\Delta^1_1(\alpha)$ in a countable ordinal $\alpha$ (note the use
of lightface font as opposed to the boldface in the previous
paragraph). We need to show that $B$ is actually a Borel set in
a countable ordinal. Now $B$ is a projection of a set in P, which
means it is the projection of a set that is $\Delta^1_1(\alpha)$
in a countable ordinal $\alpha$. Further, we have shown that $B$
is Borel. Thus $B$ satisfies the hypothesis of
Lemma~\ref{lemma:proj_is_countable}. It follows that $B$ itself is
$\Delta^1_1(\beta)$ in a countable ordinal $\beta$. In other
words, $B$ is in P. Since this argument holds for any element of NP $\cap$ co-NP,  this completes the proof.
\end{proof}

\begin{theorem}
{\rm NP $\neq$ co-NP.}
\end{theorem}
\begin{proof} This follows immediately from P $\neq$ NP and P = NP $\cap$
co-NP. An alternative proof is obtained by observing that in
\cite[Section 1]{Sch} and using the notation there, the lightface
analytic set $\Delta$ which is the projection of a lightface
$G_\delta$ set ${\mathcal G}$ is in NP but not in co-NP.
\end{proof}

\begin{remark} There is an interesting interplay between the answers to the
questions  P $\stackrel{?}{=}$ NP, P $\stackrel{?}{=}$ NP $\cap$
co-NP, and NP $\stackrel{?}{=}$ co-NP. Without the result P = NP
$\cap$ co-NP, NP $\neq$ co-NP is a stronger statement than P
$\neq$ NP because it implies P $\neq$ NP, however the converse is
not true. On the other hand, if P = NP $\cap$ co-NP, then NP
$\neq$ co-NP if and only if P $\neq$ NP.
\end{remark}

\subsection{The classes P$^+$ and NP$^+$ in infinite time Turing machines}

If we regard an input $x \in\, ^{\omega}2$ as having length
$\omega_1^x$, which is the least $x$-admissible ordinal greater
than $\omega$, then the versions of the classes P and NP obtained
are labelled P$^+$ and NP$^+$. We recall the following
definitions from \cite{Sch}

\begin{definition}
Let $A \subset\, ^{\omega}2$. We say that $A$ is in ${\rm P}^+$ if there exists a Turing machine $T$ such that\\
(a) $x \in A$ if and only if $T$ accepts $x$ \\
(b) $T$ halts on all inputs after $< \omega_1^x$ many steps.
\end{definition}

\begin{definition} Let $A \subset\, ^{\omega}2$. We say that $A$ is in ${\rm NP}^+$ if there exists a Turing machine $T$ such that \\
(a) $x \in A$ if and only if $(\exists y$ such that $T$ accepts $x \oplus y)$ \\
(b) $T$ halts on all inputs $x \oplus y$ after $< \omega_1^x$ many steps.
\end{definition}

\begin{theorem}{\rm \cite[Theorem 2.13]{Sch}} ${\rm P}^+ = P_{\omega_1^{CK}} = \Delta_1^1 $
\end{theorem}

\begin{corollary}{\rm \cite[Corollary 2.14]{Sch}}
${\rm P}^+ \neq {\rm NP}^+$.
\end{corollary}

We now show that the class ${\rm P}^+$ is a strict subset of ${\rm
NP}^+ \cap {\rm co-}{\rm NP}^+$.
\begin{theorem}
${\rm P}^+ \subset {\rm NP}^+ \cap {\rm co-}{\rm NP}^+$.
\end{theorem}
\begin{proof}Follows from Corollary~\ref{corollary:notless}.
\end{proof}

\begin{theorem}
${\rm NP}^+ \neq {\rm co-}{\rm NP}^+$
\end{theorem}
\begin{proof} Follows immediately from the proof of \cite[Theorem 2.13]{Sch}.
\end{proof}

\begin{acknowledgment}
I would like to thank Ralf Schindler for suggesting that I write
this paper and for his correspondence. I thank David Marker for his help
with Lemmas~\ref{lemma:proj_is_countable} and \ref{lemma:notless}. I thank Alexander Kechris for providing a proof sketch for Lemma~\ref{lemma:proj_is_countable}.
I also thank Albert Visser for several interesting conversations.
 \end{acknowledgment}

\end{article}

\begin{references}

\bibitem{HamLew} J. Hamkins and A. Lewis, ``Infinite time Turing machines," J. Symb. Logic, {\bf 65} (2000), 567-604.

\bibitem{Sch} R. Schindler, ``P $\neq$ NP for infinite time Turing machines," Monatshefte fur Mathematik (to appear).


\bibitem{Sip1} M. Sipser, ``The History and status of the P versus NP question," $24^{th}$ Annual ACM Symposium on theory of computing, ACM press, 603-618.

\bibitem{Mos} Y. Moschovakis, Descriptive set theory (Studies in Logic and the foundations of Mathematics, {\bf 100}), North-Holland, Amsterdam, 1980.

\bibitem{Kec} A. Kechris, Classical descriptive set theory (Graduate texts in Mathematics), Springer-Verlag, New York, 1995.

\bibitem{KMS} A. Kechris, D. Marker, and R. Sami, ``$\Pi_1^1$
Borel sets," J. Symb. Logic, {\bf 54} (1989), 915-920.



\end{references}
\end{document}